\documentclass[12pt, a4paper,twoside,openright]{article}
\usepackage{latexsym}
\usepackage{indentfirst}
\usepackage{graphicx,psfrag,epsf}
\usepackage{amsmath,amssymb,amsfonts,amsthm,amstext,amscd}
\usepackage{mathrsfs}
\usepackage{fancyhdr}
\usepackage{times}
\usepackage{enumerate}
\usepackage{url} 
\usepackage{tikz}
\usepackage{pgf}
\usepackage{setspace}

\def\argmax{\mathop{\rm argmax}} 
\def\argmin{\mathop{\rm argmin}} 

\newtheorem{defi}{Definition}[section]
\newtheorem{exam}{Example}[section]

\newtheorem{prop}{Proposition}[section]

\fancypagestyle{plain}{%
\fancyhf{}%
\fancyhead[LO]{RATIO MATHEMATICA 30 (2016) pp. 3-21}
\fancyhead[RO]{ISSN:1592-7415}
\fancyfoot[C]{\thepage }
}

\begin{document}
\title{\textbf{A geometric view on Pearson's correlation coefficient and a generalization of it to non-linear dependencies}}
\author{Priyantha Wijayatunga\\
\small Department of Statistics, Ume\r{a} School of Business and Economics, \\ \small Ume\r{a} University, Ume\r{a} 901 87, Sweden\\
\small priyantha.wijayatunga@umu.se}

\date{}
\maketitle

\begin{abstract}
Measuring strength or degree of statistical dependence between two random variables is a common problem in many domains. Pearson's correlation coefficient $\rho$ is an accurate measure of linear dependence. We show that $\rho$ is a normalized, Euclidean type distance between joint probability distribution of the two random variables and that when their independence is assumed while keeping their marginal distributions. And the normalizing constant is the geometric mean of two maximal  distances; each between the joint probability distribution when the full linear dependence is assumed while preserving respective marginal distribution and that when the independence is assumed. Usage of it  is  restricted to linear dependence because it is based on  Euclidean type distances that are generally not metrics and considered full dependence is linear. Therefore, we argue that if a suitable distance metric is used while considering all possible maximal dependences then it can measure any non-linear dependence.  But then, one must define all the full dependences.  Hellinger distance that is a metric can be used as the distance measure between probability distributions and obtain a generalization of $\rho$ for the discrete case.\\



doi: 10.23755/rm.v30i1.5;   This is the author's version of the publication
\end{abstract}

\section{Introduction}
\label{sec:intro}

Measuring association between two random quantities is of interest in many types statistical analyses and applications in various disciplines. Pearson's product moment correlation coefficient is the standard in statistical textbooks and applications for measuring linear association. And Spearman's rank correlation coefficient is capable of measuring any monotonic dependence between two random variables.  For two ordinal variables  Cram\'{e}r's V-statistic is widely used whereas Tchuprow's T-statistic is less-known and therefore less often used (see \cite{BW13} and references therein).  Furthermore, there are many other kinds of  dependence measures used  in  statistical literature, especially in applied statistical analyses. In statistical genetics for evaluation of linkage disequilibrium between genetic markers, authors of \cite{SC02} use volume tests that are discussed in \cite{DE85} as  a measures of dependence between ordinal variables with fixed margins. For massive datasets in \cite{SB2013} it is used mutual information dimension that is defined in terms of information dimension descried in \cite{RA1970}.

In \cite{BL09} it is said that \emph{``although it is customary in bivariate data analysis to compute a correlation measure of some sort, one number (or index) alone can never fully reveal the nature of dependence; hence a variety of measures are needed"}. It is also stated therein that \emph{``if (two quantities are) not totally dependent, then it may be helpful to find some quantities that can measure the strength or degree of dependence between them"}. In this article we try to develop a measure that can indicate `the' degree or strength of association between two discrete variables. Our measure can be seen as a generalization of the Pearson's correlation coefficient $\rho$ using a suitable distance metric between joint probability distributions, instead of simple Euclidean type distances that are used in $\rho$ (see below).  Given the joint probability distribution (jpd) of two discrete variables, say, $X$ and $Y$, the degree of dependence (also called association) between them is expressed as the normalized distance between  the jpd of them and that of when the independence of them is assumed. The associated normalizing constant is geometric mean of  distances between the latter and all possible jpds where full dependence between $X$ and $Y$ is assumed while retaining each marginal distribution at a time.  These latter distances are in fact the maximal distances since we obtain them by assuming full dependence. In the following we show that the Pearson's correlation coefficient is measure of this nature based on some Euclidean type distances. That is, it is the ratio of the distance between dependence and independence, and  the geometric mean of the distances that are between full linear dependences and independence. Therefore, our measure can be regarded as a generalization of $\rho$ using a suitable distance between probability distributions and considering non-linear dependencies. One thing that $\rho$ shows us is that if we need to define a strength of a dependence then we must  find or hypothesize the full dependence(s) corresponding to the given dependence.  This aspect can make numerical evaluation of the measure algorithmic or computational  since sometimes it may not be possible to obtain the full dependences easily. However, here we do not deal with such computational issues but our consideration is on defining a measure following the structure of $\rho$. For a given dependence (in terms of a jpd) finding efficiently  related jpds representing the full dependences that preserve either of the marginal is an open problem.

 First we show that, in the simple case of binary $X$ and $Y$, the $\rho$ measures the degree of dependence with a certain type of Euclidean distance, but for multinary case (and also for continuous variables) a distance in terms another type of Euclidean area is used.  But these Euclidean type distances are appropriate for measuring only linear dependences. Since we are interested in measuring any non-linear dependence we propose to use  Hellinger distance between joint probability distributions, that is called as Matsusita distance in the discrete (see \cite{MK55}). The Hellinger distance is a metric and it possesses the so-called linear invariance properties, so it is more suitable for measuring  distances between the probability distributions. Therefore, it can be used to measure any type of dependence. 

\section{Pearson's correlation coefficient $\rho$}

For random variables $X$ and $Y,$ the Pearson's correlation coefficient $ \rho(X,Y)$ is such that $\vert \rho(X,Y) \vert \leq 1$. The equality holds if  and only if $X$ and $Y$ are fully linearly dependent and $\rho(X,Y)=0$ if they are linearly independent. And the converse of the latter  is not always  true unless $X$ and $Y$ are binary. Note that the full dependence is linear in the binary (also called $2 \times 2$) case where then the $\rho(X,Y)$ is often called $\phi$-coefficient. 
 
\subsection{$2 \times 2$ case: $\phi$-coefficient}

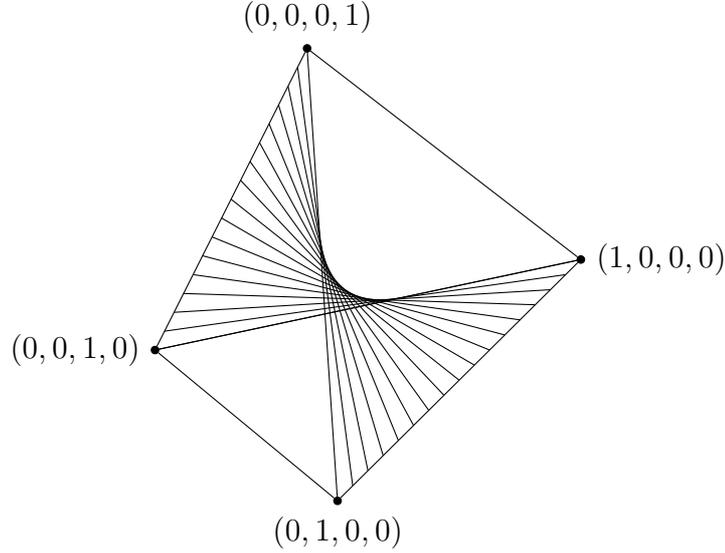
\begin{figure}
\begin{center}
\begin{tikzpicture}
[scale=4, vertices/.style={draw, fill=black, circle, inner sep=1pt}]

\node[vertices, label=below:{$(0,1,0,0)$}] (a) at (0,0) {};
\node[vertices, label=right:{$(1,0,0,0)$}] (b) at (0.8,0.8){};
\node[vertices, label=left:{$(0,0,1,0)$}] (c) at (-0.6,0.5) {};
\node[vertices, label=above:{$(0,0,0,1)$}] (d) at (-0.1,1.5) {};

\foreach \to/\from in {a/b,a/c,a/d,b/d,b/c,c/d}
\draw [-] (\to)--(\from);

\draw [-] (0.05,0.05)--(-0.13125,1.4375);
\draw [-] (0.10,0.10)--(-0.1625,1.375);
\draw [-] (0.15,0.15)--(-0.19375,1.3125);
\draw [-] (0.20,0.20)--(-0.225,1.25);
\draw [-] (0.25,0.25)--(-0.25625,1.1875);
\draw [-] (0.30,0.30)--(-0.2875,1.125);
\draw [-] (0.35,0.35)--(-0.31875,1.0625);
\draw [-] (0.40,0.40)--(-0.35,1);
\draw [-] (0.45,0.45)--(-0.38125,0.9375);
\draw [-] (0.50,0.50)--(-0.4125,0.875);
\draw [-] (0.55,0.55)--(-0.44375,0.8125);
\draw [-] (0.60,0.60)--(-0.475,0.75);
\draw [-] (0.65,0.65)--(-0.50625,0.6875);
\draw [-] (0.70,0.70)--(-0.5375,0.625);
\draw [-] (0.75,0.75)--(-0.56875,0.5625);
\draw [-] (0.80,0.80)--(-0.6,0.5);

\end{tikzpicture}
\end{center} 
\caption{Probability simplex for binary $X$ and $Y$ where their jpd  $P=(p_{00},p_{10},p_{01},p_{11})$  is a point in it. Any jpd  on surface shown by lines  represents independence of $X$ and $Y.$ \label{fig:probsimplex}}
\end{figure}

Let $X$ and $Y$ be two binary variables with a common state space $\{0,1\}$ where their jpds and marginal probability distributions are written as $p_{xy}=p(X=x,Y=y),$ $p_x=p(X=x)$ and $q_y=p(Y=y)$ for $x,y=0,1$. Let 
$P=\left(%
\begin{array}{cc}
p_{00} & p_{01} \\
p_{10} & p_{11}   \\
\end{array}%
\right)
$
for short.
As shown in \cite{FG1970}, any such $P$ can be represneted as a point in the probability simplex shown in the Figure \ref{fig:probsimplex}. The jpd of $X$ and $Y$ under the assumption that they are independent while keeping the marginal distributions fixed is 
$P^I=\left(%
\begin{array}{cc}
p_0q_0 & p_0q_1 \\
p_1q_0& p_1q_1 \\
\end{array}%
\right)
$ and the set of such probability distributions for all $P$ makes a surface (shown by lines) in the probability simplex. The $\phi$-coefficient of $X$ and $Y$ is defined by 
\begin{equation*}
\phi = \frac{p_{11}-p_1q_1}{\sqrt{p_1(1-p_1)q_1(1-q_1)}},
\end{equation*}
which is a measure of degree of association between $X$ and $Y$. Now let $X$ and $Y$ be positively correlated, then there are two jpds under the assumption that the two variables are fully dependent.  They are 
$
P^X=\left(%
\begin{array}{cc}
p_0 & 0 \\
0 & p_1 \\
\end{array}%
\right)
$
and  
$
P^Y=\left(%
\begin{array}{cc}
q_0 & 0 \\
0 & q_1
\end{array}%
\right)
$,
where $P^X$ is when the marginal distribution of $X$  is preserved and $P^Y$ is when the marginal distribution of $Y$ is preserved. Note that each full dependence is obtained from $P$ while preserving respective marginal distribution, then the marginal distribution of the other variable should be assumed by it. Therefore in these cases, the full dependence is essentially linear.

For a generalization of $\rho$ to measure `any' type of dependence we need to  look at its structure and construction. First we consider the case of two binary variables by examining  the $\phi$-coefficient. Let $D_{P^I,P}$ be $p_{11}-p_1q_1$  that is the  $(2,2)^{th}$ component Euclidean distance between the two probability distributions $P^I$ and $P$. It is a measure of how far the dependence (under $P$) from the independence (under $P^I$) when marginals of $X$ and $Y$ are fixed. Note that in the $2 \times 2$ case it is sufficient to consider a single component difference (between the two  probability matrices) since all the components have same absolute difference. Similarly, we have $D_{P^I,P^X}=p_1(1-q_1)$ and $D_{P^I,P^Y}=q_1(1-p_1)$. Since $P^X$ and $P^Y$  are the two full dependences that we can obtain from $P$ while preserving respective  marginal in each case, we have that $D_{P^I,P}\leq D_{P^I,P^X} $ and $ D_{P^I,P} \leq D_{P^I,P^Y}.$   In fact $D_{P^I,P} = p_{11}-p_1q_1=p_1(p_{11}/p_1-q_1)\leq p_1(1-q_1)=D_{P^I,P^X}$ since $p_1 \geq p_{11}$ and similarly the other inequality. It is easy to see that the denominator of the $\phi$-coefficient is the geometric mean of  $D_{P^I,P^X}$ and $D_{P^I,P^Y}$ (the two maximal distances) and the numerator is $D_{P^I,P}$. Therefore,  the $\phi$-coefficient can be thought of as  the normalized distance  between $P$ and $P^I$ where the normalizing constant is the geometric mean of the two maximal distances. Hence the $\phi$-coefficient is $1$ if and only if $P=P^X =P^Y$ (full dependence) and it is $0$ if and only if $P^I=P$ (independence). 

\subsection{$n \times m$ case}

Let $X$ and $Y$ be two multinary random variables where their state  spaces are $\{0,1,..,n-1\}$ and $\{0,1,..,m-1\}$ respectively for $n,m > 2$. For any given jpd of $X$ and $Y,$ $P=(p_{00},...,p_{0(m-1)}; p_{10},...,p_{1(m-1)};...;p_{(n-1)1},...,p_{(n-1) (m-1)})$ where $p_{ij}=p(X=i,Y=j)$ for $i=0,..,n-1$ and $j=1,...m-1$, we define the probability simplex, $\Delta =\{ P = (p_{ij})_{n \times m} : \sum_{ij} p_{ij}=1,p_{ij} \geq 0; i=0,1,..,n-1; j=0,1,...,m-1 \}$ similar to the case of two binary random variables. But here visualization of it is more difficult.  Recall that $\rho(X,Y)=cov(X,Y)/\sqrt{var(X)Var(Y)}$, where $cov(X,Y)=\sum_{x,y}xyp(x,y)-\sum_x xp(x) \sum_y yp(y)$ and $var(X)=\sum_{x}x^2p(x)-\{\sum_x xp(x) \}^2$. In the following we try to visualize the $\rho$ and its structure for understanding how it measures the dependence.

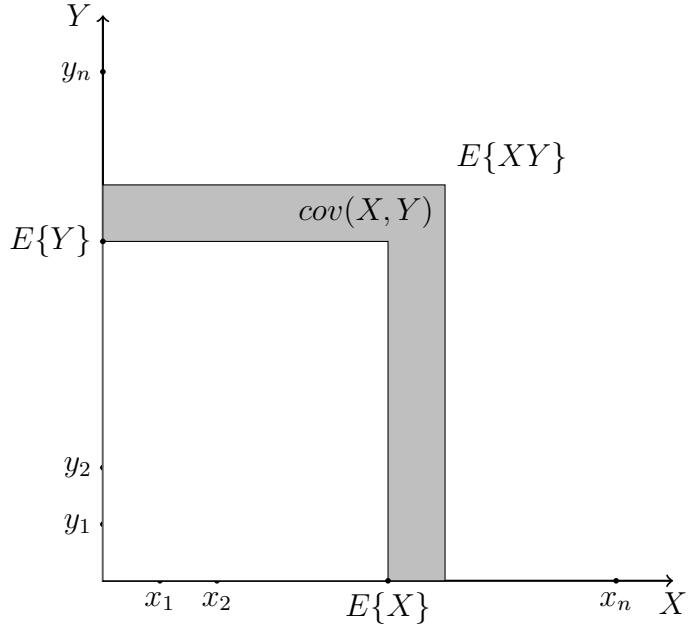
\begin{figure}
\begin{center}
\begin{tikzpicture}
[scale=3, vertices/.style={draw, fill=black, circle, inner sep=1pt}]

\draw [thick, <->] (0,2.5) -- (0,0) -- (2.5,0);
\node [below] at (2.5,0) {$X$};
\node [left] at (0,2.5) {$Y$};

\draw[fill] (0.25,0) circle [radius=0.01]; \node [below] at (0.25,0) {$x_1$};
\draw[fill] (0.5,0) circle [radius=0.01]; \node [below] at (0.5,0) {$x_2$};
\draw[fill] (2.25,0) circle [radius=0.01]; \node [below] at (2.25,0) {$x_n$};

\draw[fill] (0,0.25) circle [radius=0.01]; \node [left] at (0,0.25) {$y_1$};
\draw[fill] (0,0.5) circle [radius=0.01]; \node [left] at (0,0.5) {$y_2$};
\draw[fill] (0,2.25) circle [radius=0.01]; \node [left] at (0,2.25) {$y_n$};

\draw [fill=lightgray] (0,0) rectangle (1.5,1.75);
\node[above right] at (1.5,1.75) {$E\{XY\}$};

\draw [fill=white] (0,0) rectangle (1.25,1.5);
\draw[fill] (1.25,0) circle [radius=0.01]; \node[below] at (1.25,0) {$E\{X\}$};
\draw[fill] (0,1.5) circle [radius=0.01]; \node[left] at (0,1.5) {$E\{Y\}$};

\node[below left] at (1.5,1.75) {$cov(X,Y)$};

\end{tikzpicture}
\end{center} 
\caption{Covariance of $X$ and $Y$ is the weighted averaged Euclidean area difference. \label{fig:second}}
\end{figure}

Let us take the case where $n=m$, thus allowing us to have perfect (one-to-one) dependence between $X$ and $Y,$ linear or non-linear.  It can be seen that when $X$ and $Y$ are assigned to two perpendicular axes, $cov(X,Y)$ is area difference between two rectangular  Euclidean areas, that is shown as the dark area  in the Figure \ref{fig:second}. The first area (i.e., $ \sum_{x,y}xyp(x,y)$) is the weighted average area created by the values of  $X$ and $Y,$ where, for each component area that is being weighted is with side lengths $X=x$ and $Y=y$ and its weight is the respective joint probability of $X=x$ and $Y=y$, i.e., $p(X=x,Y=y)$. This area represents the dependence between $X$ and $Y$. And the second area  (i.e., $\sum_{x}xp(x)  \times\sum_{y}yp(y)$) is the area created by the side lengths that are the weighted average of values of $X$ (i.e., $E\{X\}$) and that of $Y$ (i.e., $E\{Y\}$) where  the weights are the respective marginal probabilities. Since the lengths or values  $E\{X\}$ and $E\{Y\}$ are also on same axes as $X$ and $Y$ are, respectively, we can see the difference of the two areas. Note that it can be seen that the second area (i.e., $\sum_{x,y}xyp(x)p(y)$) is also calculated in the similar way as  the first, but assuming the independence of $X$ and $Y$, i.e., it is the weighted average area created by the values of  $X$ and $Y$, where for each component  area that is being weighted is with side lengths $X=x$ and $Y=y$  and the weight associated with it is the respective joint probability of $X=x$ and $Y=y$ assuming independence $p(X=x,Y=y)=p(X=x)p(Y=y)$. So the second area represents the scenario of the independence of $X$ and $Y$. Therefore one can view  that the two areas refer to those  when a dependence between $X$ and $Y$ is assumed and when their independence is assumed while keeping the marginal distributions fixed, therefore $cov(X,Y)$ is a `distance' in terms of a Euclidean area difference between dependence and independence of the two variables.   

Moreover $var(X)$ can be interpreted in the same way. Now $X$ is assumed to be on both axes meaning that $Y$ is replaced by $X$ (taken as if $Y$ were $X$). This is a context of assuming a full dependence of $X$ and $Y$ when the marginal of $X$ is preserved. Assuming one  variable by the other is `a way' to consider a case of full dependence between the two variables. Then we are  assuming the marginal of $Y$ by that of $X$. This assumption is easily seen when both variables have same sizes in their state spaces but it is hard to see when they are different. So the $E\{X^2\}$ is indicated by the weighted average area that we obtain when $Y$ is $X$ where weight for each component area $x^2$ is $p(x,y)=p(x)$, i.e., when the marginal of $X$ is preserved.  This is a sensible area under full dependence. And $E\{X\}^2$ is indicated by  the area when the respective weight is $p(x)p(y)=p(x)^2$ where $x=y$. This is a hypothetical case where it is taken as if  $Y$ were $X$, yet their joint probability is taken as if  they were independent.   So, $var(X)$ is deviation of the full dependence from independence if $Y$ were $X$. And the same interpretation applies for $var(Y)$.

Thus, $\rho(X,Y)$ is the normalized area difference referring to $cov(X,Y)$ with the normalizing constant being the geometric mean of the two maximal area differences referring to $cov(X,Y)$ where they are such that, one is when $Y$ is assumed to be $X$ (i.e., $var(X)$) and the other is when $X$ is assumed to be $Y$ (i.e., $var(Y)$). That is, the normalizing constant is obtained by assuming the full dependence between $X$ and $Y.$ However the full dependence quantified  in this way is appropriate only 
for doing so for linear dependences. Since there are two such cases of full linear dependence the geometric mean of these two maximal area differences  is taken. Note that the above interpretation is valid for the case of $X$ and $Y$ have continuous state spaces. 

One thing that we need to show is that $cov(X,Y)$ is maximal (or minimal) when $X$ and $Y$ are  strictly monotonically related, for example, linearly related positively (negatively), among all cases of full ono-to-one dependencies between $X$ and $Y$ for fixed maginals of $X$ and $Y$. This indicates that $\rho$ is not able to identify non-monotonic relations since their covariance values can not be ordered. To see that $cov(X,Y) $ is maximal when $Y$ is strictly increasing with $X$, let  $\mathcal{X}=\{a_1<...<a_n\}$ be the state space of $X$ and $\mathcal{Y}=\{b_1<...<b_n\}$ be that of $Y$. Then considering inequalities $(a_i-a_j)(b_i-b_j) >0$  for $i,j=1,...,n$ (i.e., we have $a_ib_i+a_jb_j > a_ib_j +a_jb_i$) it can be shown that $\sum_i a_ib_i > \sum_{i,j:j=f(i)}a_ib_j $ where $f$ is any one-to-one function from $\mathcal{X}$ to $\mathcal{Y}$ such that $f(i) \neq i$  for at least two distinct values of  $i$ (i.e., $f$ is not a strictly increasing function of $i$). Now if the marginals of $X$ and that of $Y$ are $(p_1,...,p_n)$ and $(q_1,...,q_n)$, where $p_i=q_i$ for all $i=1,...,n$ when $Y$ is monotonically increasing with $X$ and otherwise $p_i=q_j$ for some appropriate $i \neq j$ for $i,j=1,....,n$, then  $\sum_i a_ib_i p_i > \sum_{i,j:j=f(i)}a_ib_j p_i$ meaning that $E\{XY_M\} \geq E\{XY\}$ where $Y_M$ is $Y$ when it is strictly increasing with $X$. This implies that $cov(X,Y_M) \geq cov(X,Y)$ for fixed marginals of $X$ and $Y$. Therefore, for discrete $X$ and $Y$, $ \rho(X,Y) $ is maximal when $Y$ is strictly increasing in $X$, among all one-to-one relationships between them. So, if this is the case $ \rho(X,Y) =1$ (maximal) since $cov(X,Y) \leq var(X) $  and $cov(X,Y) \leq var(Y).$

\section{Some other popular measures of dependence}

There are a few  popular measures of dependence  that have similar structure in their definition. We review them briefly by giving some interpretations that support our definition of dependence measure.   

\subsection{Spearman's rank correlation coefficient $\rho^s$}

In many statistical analyses, especially for non-normal data a popular measure of dependence between two random variables, say, $X$ and $Y$, is the Spearman's rank correlation coefficient. 

\begin{equation*}
\rho^s = 1 - \frac{6 \sum_{i=1}^n d_i^2}{n(n^2-1)}
\end{equation*}
where $d_i = x_{(i)} - y_{(i)}$  and $x_{(i)}$ is the $i^{th}$ smallest value in the data sample of $X$ and similarly for $y_{(i)}$. It is obvious that $\rho^s = 1$ if and only if two components of  data pair $(x_i,y_i)$ has the same ranking, for all data pairs since then $d_i=0$ for all $i$. And one can see that for a perfect negative dependence $\sum_{i=1}^n d_i^2 $ should be its maximal value that is $n(n^2-1)/3$ in order to get     $\rho^s_{X,Y}= -1$.  Therefore the  normalizing constant is taken as  $n(n^2-1)/6$ but due to the structure of the definition of the coefficient it is applied to the term $\sum_{i=1}^n d_i^2$.  Therefore the  $\rho^s$ is an accurate   measure any monotonic dependence between the two variables. However, when the two variables are not having a strictly monotonic relationship the measure can not give a correct picture of the dependence. 

\subsection{Information theoretic measures}

Another popular measure of dependence, especially in machine learning literature and applied statistics is so-called mutual information (see, for example, \cite{WP2006}). For discrete random variables $X$ and $Y$, it is defined as
\begin{equation*}
I(X,Y)= \sum_{x,y} p(x,y)log \frac{p(x,y)}{p(x)p(y)}
\end{equation*}
and furthermore, conditional mutual information between $X$ and $Y$ given another variable $Z$ is defined as
\begin{equation}
CI(X,Y,Z)= \sum_{x,y,z} p(x,y,z)log \frac{p(x,y \vert z)}{p(x \vert z)p(y \vert z)}
\end{equation}
If  $X$ and $Y$ are independent then the $I(X,Y)=0$ and if $X$ and $Y$ are conditionally independent given $Z$ then the $CI(X,Y,Z)=0$.  In fact, these dependence measures are also based on so-called Kullback-Leibler (KL) distance or rather divergance, \cite{KL1951}. It is easy to see that $I(X,Y)$ is the KL divergence between the joint probability distribution of $X$ and $Y$, and that when independence is assumed, therefore it measures the dependence in terms of `departure' from independence.  In fact, $I(X,Y)$ is the weighted average of Euclidean distance between logarithmic of the joint probability $p(x,y)$ and that when independence is assumed, where weights are the respective joint probabilities.  That is, it is the expectation, under the joint probability, of the difference between the logarithmic of the joint probability $p(x,y)$ and that when independence is assumed. Note that though $0 \leq I(.,.) \leq 1$, there is no normalization (with respect to any maximal dependence) is involved. 

Though these information measures are used to identify respective dependences they are not  metrics since KL-divergance is not a true distance (metric), therefore they can not be used to measure the degree of dependence between variables. For example, as shown in \cite{SV1998} let $p(x,y)$ and $q(x,y)$ define two dependencies between $X$ and $Y$ where
$p(x,y)=\left(%
\begin{array}{cc}
3/8 & 1/8 \\
1/8 & 3/8 \\
\end{array}
\right)%
$
and 
$q(x,y)=\left(%
\begin{array}{cc}
1/2 & 0 \\
1/8 & 3/8 \\
\end{array}
\right)%
$.
Obviously probability distribution $q$ shows a higher dependency than that of $p$ but its mutual information is lower than that of $p$, $(MI_p(X,Y) > MI_q(X,Y))$. Note that $q$ is obtained from $p$ without preserving the marginal distributions of $X$ and $Y$. Now let $r(u,v)$ and $s(u,v)$ define two dependencies between random variables $U$ and $V$ where
$r(u,v)=\left(%
\begin{array}{ccc}
 0 & 1/7 & 1/7 \\
 1/7 & 1/7 & 1/7 \\ 
 1/7 & 1/7 & 0  \\
\end{array}
\right)%
$
and 
$s(u,v)=\left(%
\begin{array}{ccc}
 0 & 0 & 2/7 \\
 1/7 & 2/7 & 0 \\ 
 1/7 & 1/7 & 0 \\
\end{array}
\right).
$
Then we have that $MI_r(U,V) < MI_s(U,V)$. Note that $s$ shows a higher dependency than that of $r$ and it is obtained from $r$ by preserving the marginal distributions of $U$ and $V$. Furthermore, all zeros in $r$ are also in $s$. If this is the case then higher dependency implies higher mutual information.  So mutual information is restricted measure of degree of dependence.

\subsection{Chi squared test statistic $\chi^2$}

We can see that well-known Chi squared test statistic $\chi^2$ that is used for testing independence of two discrete random variables uses a certain dependence measure in it for performing the test. Let $X$ and $Y$ take values $i=1,...,\alpha$ and $j=1,...,\beta$, respectively and let us write the joint probability of $X=i$ and $Y=j$ as $p_{ij}$, marginal probability of $X=i$ as $p_{i.}$ and  that of $Y=i$ as $p_{.j}$. So, the conditional probability of $X=i$ given $Y=j$ is $p_{i \vert j}=p_{ij}/p_{.j}$ and similarly $p_{j \vert i}$ is defined. Then,
\begin{align*}
\chi^2 &=\sum_{i,j} n \frac{(p_{ij}-p_{i.}p_{.j})^2}{p_{i.}p_{.j}} = n \Big\{ \sum_{i,j} \frac{p_{ij}^2}{p_{i.}p_{.j}} -1 \Big\} =  n \Big\{ \sum_{i,j} p_{ij} \frac{p_{ij} -p_{i.}p_{.j}}{p_{i.}p_{.j}} \Big\}    \\
         & =  n \Big\{ \sum_{i,j} p_{ij} \frac{p_{i \vert j} -p_{i.}}{p_{i.}} \Big\}   = n \Big\{ \sum_{i,j} p_{ij} \frac{p_{j \vert i} -p_{.j}}{p_{.j}} \Big\}  = n E \{A\}
\end{align*}
where $A$  is a random variable taking the value $\frac{p_{i \vert j} -p_{i.}}{p_{i.}}=\frac{p_{j \vert i} -p_{.j}}{p_{.j}}$ with probability $p_{ij},$  for $i=1,...,\alpha$ and $j=1,...,\beta$, and $E$ denotes the expectation. That is, $\chi^2$ is $n$-multiple of the expectation of a random variable whose $(i,j)^{th}$ value is a `normalized' distance between the probability value $p_{i \vert j}$ and $p_{i.}$ where the normalizing constant is $p_{i.}$, for all $i,j,$  and vice versa. Note that $\frac{p_{i \vert j} -p_{i.}}{p_{i.}}$ may be referred to as the `degree' of dependence between the two events $X=i$ and $Y=j$.  In fact, it is the certainty factor for the case $p_{i \vert j} < p_{i.},$ as described in \cite{BF2002} for measuring the dependency between the two events and it is a symmetric measure. However, here it is used without the condition. So, $E\{A\}$ is the expectation of a degree of dependence between the events $X=x$ and  $Y=y$ for all $x,y.$  Therefore, $E\{A\}$ can be thought of as measure of degree of dependence between $X$ and $Y.$ And the term $n$ in $\chi^2$ makes it a statistic. That is, a statistic for testing dependence between two variables can be seen as a product of two factors; one is a quantity related the degree of dependence between two variables and the other is that of total number of data cases that are used to estimate the probabilities related to them (i.e., sample information).

\subsection{Test of two proportions}

Sometimes one may be interested in testing equality of two proportions to see if given  two variables are independent, for example, when the outcome ($Y$) of interest is binary, such as voting, denoted by $Y=1$ (or not, denoted by $Y=0$), for a political candidate in an election for two groups/populations ($X$) such as men, denoted by $X=1,$ and women, denoted by $X=0$. Then one can test if two  proportions are equal, i.e., $p(Y=1 \vert X=1)=p(Y=1 \vert X=0)$ (let us write it as $p_1=q_1$) by the $Z$ statistics 

\begin{equation*}
Z= \frac{1}{\sqrt{1/a+1/b}} \frac{p_1-q_1}{\sqrt{p(1-p)}}
\end{equation*}
where $a$ and $b$ are the sizes of the two samples of $Y$ when $X=1$ and $X=0$, respectively, and $p=p(Y=1)$.  Now we can interpret that the factor $\frac{p_1-q_1}{\sqrt{p(1-p)}}$ as a   measure of degree of dependence between the two variables due to the term $(p_1-q_1)$ in it, where the term $\sqrt{p(1-p)}$ should be taken as the normalizing constant. Note that the latter is constructed assuming full dependence between the two variables where, then their joint probability distribution is
$
P=\left(%
\begin{array}{cc}
1-p & 0 \\
0 & p \\
\end{array}
\right)%
$
or similar. Instead of just using $p$ which is the pooled proportion,  the geometric mean of $p$ and $(1-p)$ should be used as the normalizing constant. This is necessary to yield the same test statistic value for testing the same hypothesis with complementary probabilities i.e., $p(Y=0 \vert X=1)$ and $p(Y=0 \vert X=0)$. And the term $\frac{1}{\sqrt{1/a+1/b}} $ which is a function of sample sizes  (sample information) makes $Z$ a statistics. So, similar to $\chi^2$ statistic, $Z$ has a measure of degree of dependence between the two variables in it, in addition to information on the sample sizes. 

\section{Axioms of an ideal measure of dependence}

Before we define our measure of strength/degree of dependence (or rather a generalization of $\rho$) it is appropriate to mention axioms that an ideal measure should possess as shown in \cite{GC2004}.  However, it is hard to find dependence measures satisfying all these axioms. Our generalization of $\rho$ seems to have a bigger potential in satisfying them, but we omit the discussion here. Following are the axioms;

\begin{enumerate}
\item It is well-defined for both continuous and discrete case
\item It is normalized such that its value $0$ implies the independence and value $1$ implies the full dependence (one variable is a deterministic function of the other), where all intermediate degrees of dependencies lie between $0$ and $1$
\item It is equal or has a simple relationship with the Pearson's correlation coefficient in the case of a bivariate normal distribution
\item It is a metric, i.e., it is a true measure of distance (between the independence and dependence of interest) not just a divergence
\item It is invariant under continuous and strictly increasing transformations. 
\end{enumerate}
These axioms are straightforward and require no further explanation. 

In the following we define our measure following the structure and the construction of $\rho$ but using a true distance metric. We propose to use so-called  Hellinger distance but one may use another suitable distance metric. Since we are keeping the structure of the $\rho$ the same but replacing its distance measure with a better one (a metric) when defining our dependence measure, we call it as a generalization of the $\rho$. This means that for any given dependence we should be able to define the corresponding all possible full dependences, since the measure should be a ratio between a distance from independence to the given dependence and  geometric average of distances from independence  to  the full dependences.

\section{Defining a measure of degree of dependence}

As we have seen earlier, in the two binary variables $(2 \times 2)$ case where only the linear dependence exists the dependence can be measured by using a single component Euclidean distance between joint probability distributions. However, in the case of two multinary variables ($n \times n$, where $n>2$) we can have many types of dependences, and therefore distances among probability distributions can not be defined through only a single component or  a weighted average area difference, that are Euclidean type distances and capable of measuring only linear dependences.  Therefore we need to use some other suitable distance to measure any non-linear dependences.   In the following we discuss a possible distance that is a true metric.

\subsection{A metric distance between two probability distributions}

We propose to use Hellinger distance between probability distributions (also called Matsushita distance for the discrete case) which is a metric in the probability simplex for our task of measuring dependence. Recall that our dependence measure should be the normalized distance between the given joint probability distribution of the two variables and that when their independence is assumed while preserving the marginals, where the normalizing constant is obtained by considering similar distances related to the all possible maximal dependences but preserving only one of the marginals at each time.  Let  $\Phi$ and $\Psi$ be two discrete distribution functions ($\phi$ and $\psi$ are probability distributions or mass functions) then the Hellinger distance between $\Phi$ and $\Psi$ is defined as
\begin{equation*}
M(\Phi,\Psi) = \bigg\{ \frac{1}{2}\sum_{x} \bigg\{  \sqrt{\phi(x)} - \sqrt{\psi(x)} \bigg\}^{2} \bigg\}^{1/2}
\end{equation*}
In addition to satisfying properties of a metric $M(.,.)$ also satisfies the following properties: (1) $0 \leq M(\Phi,\Psi) \leq 1$, (2) $M(\Phi(T),\Psi(T)) =M(\Phi(T+a),\Psi(T+a)) $ for any constant $a$, and (3) $M(\Phi(T),\Psi(T)) =M(\Phi(cT),\Psi(cT)) $ for any constant $c \neq 0$ where the last two are called the \emph{linear invariance} properties of the probability metric. Note that $\big( M(.,.) \big)^2$ is not a metric. 

First we should have an idea about the furtherest jpd(s) for a given jpd that may represent independence. In fact we can see that the furtherest probability distribution to a distribution that represent independence is not useful but those with fixed marginals, each at a time. For a given distribution function, say, $\Phi$ let us find the maximally Hellinger-distanced  distribution function $\Psi$. The following proposition shows how to find it. 
\begin{prop} 
For positive probability distribution $\phi$ maximally Hellinger-distanced probability distribution $\psi$ is given by
\begin{equation*}
\psi(t)= \begin{cases} 1, & \mbox{if } t= \argmin_u \phi(u) \\ 0, & \mbox{otherwise. } \end{cases}
\end{equation*}
and then, $M(\Phi,\Psi) =\Big\{ 1- \sqrt{min \large\{ \phi(t): t \in \mathcal{T} \large\}}  \Big\}^{1/2} < 1$. 
\end{prop}
\textbf{Proof. }Let $\vert \mathcal{T} \vert = n$, $\phi(t_i)=\phi_i$ and $\psi(t_i)=\psi_i$ for $i=1,...,n$. Let re-index all $\phi_i$'s 
such that $\phi_{(1)} \geq \phi_{(2)} \geq .... \geq \phi_{(n)}$ and possibly  some of the $\psi_i$'s can be zeros. $M(\Phi,\Psi)$ is maximal when $\sum_{t \in \mathcal{T}}  \sqrt{\phi(t)\psi(t)}$ is minimal.

\begin{eqnarray*}
\sum_{i=1}^n \sqrt{\phi_i \psi_i} &=& (\sqrt{\psi_1 }+ ... + \sqrt{\psi_n})\sqrt{\phi_{(n)}}  \\
                                                  & & + (\sqrt{\psi_1} + ... + \sqrt{\psi_{n-1}})(\sqrt{\phi_{(n-1)}}-\sqrt{\phi_{(n)}}) \\
                                                 &  & ... + \sqrt{\psi_1} (\sqrt{\phi_{(1)}}-\sqrt{\phi_{(2)}}) \geq  \sqrt{\phi_{(n)}}
\end{eqnarray*} 
That is, $\sum_{i=1}^n \sqrt{\phi_i \psi_i}$  is minimal when $\psi_1 = ... = \psi_{n-1}=0$ and $\psi_n=1$. So we obtain the maximally Hellinger-distanced distribution function $\Psi$ and therefore $M(\Phi,\Psi)$.$\Box$

But then $T$ is deterministic variable with respect to $\Psi$! This theorem says that for any given probability distribution, bivariate discrete in our case, the maximally Hellinger-distanced probability distribution is represented by a vertex of the probability simplex. All its component  are zeros except for one place that has 1 that is corresponding to the smallest probability value of the reference probability distribution.  This is a degenerate case as far as dependence of the two variables are concerned since it represents that both variables are deterministic and having full dependence. Therefore, such a full dependence can not be used for the normalization since it does not generally preserve the marginals.

For a given jpd $P$ of $X$ and $Y,$ the dependence of them that it represents should be measured with a suitable  normalized distance between $P$ and $P^I$. It is clear from above that the normalizing constant should be the geometric mean of distances from independence to all possible full dependences where each such full dependence should be preserving either of marginals. This rule is to follow the correlation coefficient definition. Therefore, an essential step is to find the two types of probability distributions $P^X$ (jpd(s) representing full dependence when marginal of $X$ is fixed) and $P^Y$ (jpd(s) representing full dependence when marginal of $Y$ is fixed) in order to find the normalizing constant. As you will see in some cases there may be multiple candidates for each of them. Therefore we have the following definition. Note that there are some instances such as in \cite{GC2004} and \cite{ST1996} where Hellinger distance between the jpd and that of when independence is assumed is used for measuring the dependence, but in such work no normalization is done. However,  the above proposition implies that distance between any non-deterministic  jpd representing independence and  that representing a full dependence can be strictly less than $1$ for two discrete random variables, therefore normalization is necessary if one wants to have a measure that shows strength of dependence.

\begin{defi}  When $M$ is a metric in the probability simplex of two discrete random variables $X$ and $Y,$ $M$-based measure of degree of dependence between $X$ and $Y$ represented by their joint distribution function $P$ is defined as 
\begin{equation*}
\rho^M(X,Y)= \frac{M(p^I,p)} { \big\{ \prod_{p^X \in \mathcal{P}^X } M(p^I,p^X)^{1/\vert \mathcal{P}^X \vert } \prod_{p^Y \in \mathcal{P}^Y }   M(p^I,p^Y)^{1/\vert  \mathcal{P}^Y \vert}  \big\}^{1/2}} 
\end{equation*}
 where $P^I$ is the joint distribution function of $X$ and $Y$ when their independence is assumed, $\mathcal{P}_{max}^X$ denotes the set of all joint distribution functions,  each representing a maximal dependence while preserving the marginal distribution of $X$ and similarly for  $\mathcal{P}_{max}^Y$,  $\vert A \vert$ is the cardinality of the set $A$, and $M(P,Q)$ is the distance metric between two probability distributions $P$ and $Q.$
\end{defi}
Note that the denominator is the geometric mean of the maximal distances between full dependences and the independence. And we use Hellinger distance as the distance measure. Since $\rho^M$ is defined following the structure of the Pearson's correlation coefficient it can be regarded as a generalization of it for the case of discrete variables.

For linear relationships measuring the dependence is relatively easy since  both $P^X$ and $P^Y$ represent perfect linear dependence. This is when they have all  their entries zero except for those, but may not be all, in each diagonal in respective case. For example, for a positive linear relation, $P^X$ is obtained by assigning each main diagonal entry with the sum of all entries in the respective row. This assures that the marginal probability of $X$ is preserved when obtaining full dependence, and similarly for $P^Y.$ Note that positive linear relationship is selected if main diagonal entries are generally larger than the other entries in the joint probability value matrix $P$. But when we allow  non-linear relationships between $X$ and $Y$ there are no pre-specified $P^X$ and $P^Y,$ therefore multiple candidates may exist for each of them.  We argue that they should be induced from the jpd  in a similar way to the case of linear dependence. So we propose following simple rule  for obtaining $P^X$ and $P^Y.$

\begin{defi}
For each $x$, when there exists a single value $y'$ such that $y'=\argmax_y p(X=x,Y=y)$, then let $p^X(X=x,Y=y')=p(X=x)$ and $p^X(X=x,Y \neq y')=0$ to obtain $P^X$.  If there are multiple such $y'$ values then obtain multiple $P^X$, each refering to one of those $y'$ values, assuming that it is the only value where maxima exists. And similarly $P^Y$ is defined. 
\end{defi}
By this way, we get one or more jpds each representing a maximal dependence that preserves respective marginal.

\section{Examples of $n \times n$ case where $n \geq  2$} 

 Now we consider some different cases of $P$ and demonstrate how we can calculate our measure and compare its value to those of some trational measures.

\paragraph{Case 1} 

Suppose a simple case of each row and column of $P$ having a single maximal entry that is common to both its row and column. Then the other entries in the row are summed onto the maximal entry in the row for each row to yield $P^X$ and similarly $P^Y$ is obtained.  Therefore, $P^X$ and $P^Y$ are on the boundary of $\Delta$, so they are the furtherest probability distributions from $P^I$ while preserving respective marginals. Then the degree of dependence between $X$ and $Y$ is defined as (since  $\vert \mathcal{P}_{max}^X \vert =\vert \mathcal{P}_{max}^Y \vert =1$ )
\begin{equation*}
\rho^M(X,Y)= \frac{M(P^I,P)}{\sqrt{ M(P^I,P^X) M(P^I,P^Y)}}
\end{equation*}
\begin{exam}
For binary $X$ and $Y$ with 
$
P=\left(%
\begin{array}{cc}
0.3 &  0.2 \\
0.1 & 0.4 \\
\end{array} 
\right)%
$,
$\phi=0.4082$ and $\rho^M=0.2783$ (Cramer's $V$ and Tschuprow's $T$ are $0.4082$). And interchanging off-main diagonal entries but keeping the main diagonal entries as they were, i.e., having
$
P=\left(%
\begin{array}{cc}
0.3 &  0.1 \\
0.2 & 0.4 \\
\end{array}
\right)%
$,
 gives the same results for all measures.
\end{exam}

\begin{exam} Let state spaces of $X$ and $Y$ be $\{1,2,3\}$ and their joint probability
$
P=\left(%
\begin{array}{ccc}
0.05 & 0.03 & 0.20 \\
0.30 & 0.07 & 0.05 \\
0.04 & 0.20 & 0.06 \\
\end{array}
\right)%
$
that is a non-linear dependence and then
$
P^I=\left(%
\begin{array}{ccc}
0.1092 & 0.084 & 0.0868 \\
0.1638 & 0.126 & 0.1302 \\
0.1170 & 0.090 & 0.0930 \\
\end{array}
\right)%
$,
$  
P^X=\left(%
\begin{array}{ccc}
0.00  & 0.00 & 0.28 \\
0.42  & 0.00 & 0.00 \\
0.00  & 0.30 & 0.00 \\
\end{array}
\right)%
$ 
and
$
P^Y=\left(%
\begin{array}{ccc}
0.00  & 0.00 & 0.31 \\
0.39  & 0.00 & 0.00 \\
0.00 &  0.30 & 0.00 \\
\end{array}
\right)%
$. 
And then $\rho=-0.2025$ but $\rho^M=0.4113$  (Cramer's $V$ and Tschuprow's $T$ are $0.5472$). But had that 
$
P==\left(%
\begin{array}{ccc}
0.05 & 0.03 & 0.20 \\
0.04 & 0.20 & 0.05 \\
0.30 & 0.07 & 0.06  \\
\end{array}
\right)%
$
which is a linear dependence then  $\rho=-0.5474$ and $\rho^M=0.4075$  (Cramer's $V$ and Tschuprow's $T$ are $0.5467$). Note the change in the degree of dependence is small since linear dependence is obtained from nonlinear case by just interchanging probability values in $P$.
\end{exam}

\paragraph{Case 2}

When each row and column of $P$ has a single maximal entry that may not be common to both its row and column we still can obtain a single $P^X$ and a single $P^Y$. Therefore, we can apply the above definition.
\begin{exam}
When
$
P=\left(%
\begin{array}{ccc}
0.30 & 0.03 & 0.20 \\
0.05 & 0.07 & 0.05 \\
0.04 & 0.20 & 0.06 \\
\end{array}
\right)%
$
we have $\rho=0.1383$ and $\rho^M=0.450011$. Note that here we have that Cramer's $V$ and Tschuprow's $T$ are $0.4257843$  that are lesser than our measure.
\end{exam}

\paragraph{Case 3} 

When there are more than one maximal entry in a row or a column we have multiple $P^X$'s and multiple $P^Y$'s. Note that here we try to obtain a similar situation in the above two cases. That is, each row of $P^X$ has only one non-zero element (it is obtained by summing up all entries in the corresponding row of $P$, thereby preserving the marginal probability distribution of $X$).  Assume that we get $a$ number of $P^X$'s, say, $P^{X_1},...,P^{X_a}$ and $b$ number of $P^Y$, say,  $P^{Y_1},...,P^{Y_b}$.  Let us consider the following example.

\begin{exam}
When
$
P=\left(%
\begin{array}{ccccc}
0.11 & 0.01 & 0.01 & 0.01 & 0.01 \\
0.01 & 0.01 & 0.01 & 0.01 & 0.25 \\
0.01 & 0.10 & 0.10 & 0.01 & 0.01 \\
0.01 & 0.01 &  0.01 & 0.15 & 0.01 \\
0.01 & 0.10 & 0.01 & 0.01 & 0.01 \\
\end{array}
\right)%
$
then 
%
 we make two $P^X$'s;

$
P^{X_1}=\left(%
\begin{array}{ccccc}
0.15 & 0.000 & 0.000 & 0.00 & 0.00 \\
0.00 & 0.000 & 0.000 & 0.00 &  0.29 \\
0.00 & 0.230 & 0.000 & 0.00 & 0.00  \\
0.00 & 0.000 & 0.000 & 0.19 & 0.00 \\
 0.00 & 0.140 & 0.000&  0.00 & 0.00 \\
\end{array}
\right)%
$
and  \\
$
P^{X_2}=\left(%
\begin{array}{ccccc}
0.15 & 0.000 & 0.000 & 0.00 & 0.00 \\
0.00 & 0.000 & 0.000 & 0.00 &  0.29 \\
0.00 & 0.000 & 0.230 & 0.00 & 0.00  \\
0.00 & 0.000 & 0.000 & 0.19 & 0.00 \\
0.00 & 0.140 & 0.000 &  0.00 & 0.00 \\
\end{array}
\right)%
$.\\
Therefore we have two maximal distances to these two full dependences. They are  $M(P^I,P^{X_1})$ and $M(P^I,P^{X_2})$ and similarly we obtain another two full dependences when marginal of $Y$ is preserved. Therefore, 
\begin{equation*}
\rho^M(X,Y)= \frac{M(P,P^I)}{ \prod_{i=1}^2 \prod_{j=1}^2 \big\{M(P^{X_i},P^I) M(P^{Y_j},P^I)\big\}^\frac{1}{4}}
\end{equation*}
Then $\rho=-0.0491$ and $\rho^M=0.5731$. Note that here we have that Cramer's $V$ and Tschuprow's $T$ are $0.6652$. 

\end{exam}
\section{Conclusion}

We have looked at the structure and the construction of the Pearson's correlation coefficient $\rho$ in order to have a generalization of it for measuring any non-linear dependence between two random variables.   We have shown that it is simple do it geometrically for discrete variables.  It can be shown that $\rho$ is a normalized  `Euclidean' type distance between the joint probability distribution of the two random variables and that when their independence is assumed in the probability simplex of the two variables  where normalizing constant is the geometric mean of two maximal such distances; each between full linear dependence of the two variables and their independence while  preserving the marginal distribution of respective variable. So, we have shown that if we consider all possible full dependences and use an appropriate distance such as Hellinger then we can have a genaralization of $\rho$. But generally it is not easy to find all possible maximal distances, which is an open problem that may need algorithmic or computational solutions. However we have shown some examples after having defined a generalization.  

\paragraph{Acknowledgments:} Financial support for this research is from Swedish Research Council for Health, Working Life and Welfare (FORTE) and Swedish Initiative for Microdata Research in the Medical and Social Sciences (SIMSAM).

{}

\end{document}